\documentclass[11pt]{amsart}

\usepackage[T1]{fontenc}
\usepackage{lmodern}
\usepackage{amsmath,amssymb,amsthm,mathtools}
\usepackage{microtype}
\usepackage{enumitem,needspace}
\usepackage[colorlinks=true,linkcolor=blue,citecolor=blue,urlcolor=blue]{hyperref}
\hypersetup{
 pdftitle={Prime-Power Constant Terms in Bayart's Power Series Problem},
 pdfauthor={Elad Paran},
 pdfsubject={Unique factorization in formal power series rings}
}
\setlist[enumerate]{label=\textup{(\roman*)},leftmargin=2.5em,itemsep=2pt,topsep=4pt}
\newtheorem{theorem}{Theorem}[section]

\newtheorem{lemma}[theorem]{Lemma}

\theoremstyle{remark}

\newcommand{\Q}{\mathbb Q}
\newcommand{\N}{\mathbb N}
\DeclareMathOperator{\Frac}{Frac}
\DeclareMathOperator{\ord}{ord}
\newcommand{\xx}{\boldsymbol{x}}

\title[Bayart's Power Series Problem]{On Bayart's Unique Factorization Problem}
\author{Elad Paran}
\address{Department of Mathematics and Computer Science, The Open University of Israel, Ra'anana 4353701, Israel}
\email{paran@openu.ac.il}
\date{\today}

\begin{document}

\begin{abstract}
Let $A$ be a ring such that the formal power series ring $A[[t]]$ is a unique factorization domain. Bayart's problem asks whether the formal power series ring in two variables $A[[x,t]]$ must also be a UFD. Equivalently, is every irreducible element $f\in A[[x,t]]$ prime? We prove that the answer is yes whenever the constant term of $f$ is an associate of $p^m$, for any prime $p \in A$ and $1\leq m\leq3$, and for $m=4$ if $A$ contains $\Q$. We also prove the following reduction: To prove a positive answer to Bayart's problem, it suffices to prove that $A[x][[t]]$ is a UFD whenever $A[[t]]$ is.
\end{abstract}

\maketitle

\section{Introduction}

Let $A$ be a unique factorization domain (UFD). The question whether
the power series ring $A[[t]]$ is again a UFD has a long history, beginning with Krull, who suspected in 1938 \cite{Krull} that the answer is negative, which was confirmed by Samuel in 1961 \cite{Samuel}. Following this pivotal work of Samuel came many papers dealing with the question of {\it when} $A[[t]]$ is a UFD \cite{Bayart80,Bayart81,Buchsbaum,Claborn,DanilovSamuel,DanilovDGC,Lipman,Salmon}. Many partial results were obtained, yet we do not have a full answer to this question in terms of the arithmetic of the ring $A$ itself. 

One particular long-standing problem is due to Bayart: 
Suppose $A[[t]]$ is a UFD, must $A[[x,t]]$ also be a UFD \cite{Bayart73}? This question has recently been answered affirmatively for Noetherian coefficient rings \cite{Tran}, but remains open in general. An equivalent formulation \cite[Lemma~2.1]{Buchsbaum} of the problem is the following: Suppose that $A[[t]]$ is a UFD. Is every irreducible element in $A[[x,t]]$ a prime?

In this paper, we prove a partial positive result without assuming that $A$ is Noetherian

\begin{theorem}\label{thm:small}
Let $A$ be a domain such that $A[[t]]$ is a UFD, let $p\in A$ be
prime, and let $f\in A[[x,t]]$ be irreducible with constant term
associated to $p^m$. If $1\leq m\leq3$, then $f$ is prime. 
The same conclusion holds if $m=4$ and $A$ contains $\Q$.
\end{theorem}

We also give a reduction of Bayart's problem to a question about
polynomial extensions.

\begin{theorem}\label{thm:reduction}
Let $A$ be a domain. If $A[u,v][[t]]$ is a UFD, then
$A[[x,t]]$ is a UFD. In particular, suppose that for every domain $B$,
\begin{equation}\label{eq:pep}
 B[[t]]\text{ is a UFD}\quad\Longrightarrow\quad
 B[x][[t]]\text{ is a UFD}.
\end{equation}
Then $A[[t]]$ being a UFD implies that $A[[x,t]]$ is a UFD for
every domain $A$.
\end{theorem}

Section~2 contains preliminary results that are used in
the proof of Theorem~\ref{thm:small}. The proof concludes in
Section~3. Section~4 proves Theorem~\ref{thm:reduction}.

\section{Preliminaries}
\label{sec:prelim}

Throughout the paper, all rings are commutative with identity.
For a power series $f$ in finitely many variables, we write $f_0$
for its constant term. We write $f\sim g$ when $f$ and $g$ are
associated in the ring under consideration.

In this section, let $A$ be a domain such that $A[[t]]$ is a UFD,
and let $p\in A$ be prime. Then $A$ is a UFD. Indeed, a strict
ascending chain of principal ideals in $A$ remains strict in
$A[[t]]$, so $A$ is atomic. An irreducible element of $A$ remains
irreducible, and hence prime, in $A[[t]]$. It is therefore prime
in $A$. Put
\[
 V=A_{(p)},\qquad W=\widehat V,\qquad k=W/pW=\Frac(A/pA).
\]
Thus $V$ and $W$ are DVRs with uniformizer $p$, and the power series
rings over $W$ in finitely many variables are regular local UFDs.

We first record a coefficient observation that will be used
throughout the proof. Here and below, $\xx$ denotes a finite tuple
of formal variables.

\begin{lemma}\label{lem:division}
One has
\[
 W\cap A[1/p]=A,\qquad p^jW\cap A=p^jA\quad(j\geq0).
\]
If $a,b\in A[[\xx]]$, $a_0=p^m$, and $a$ divides $b$ in
$W[[\xx]]$, then $a$ divides $b$ in $A[[\xx]]$.
These assertions also hold with $W$ replaced by any complete DVR
extension $W'$ in which $p$ remains a uniformizer.
\end{lemma}

\begin{proof}
For $c\in A$, the element $c/p^j$ belongs to $W$ if and only if
the $p$-valuation of $c$ is at least $j$. Since $A$ is a UFD,
this is equivalent to $c\in p^jA$. The same argument applies
over $W'$.

The quotient $b/a$ is uniquely defined in $\Frac(A)[[\xx]]$.
Comparison of coefficients, in increasing total degree, shows that
its coefficients belong to $A[1/p]$. If the quotient also belongs
to $W[[\xx]]$ or $W'[[\xx]]$, the intersection just proved puts
every coefficient in $A$.
\end{proof}

For a domain $D$, positive integers $d,n$ with $n\geq2$, and
variables $z_1,\ldots,z_d$, let
\[
 \sigma_n:D[[z_1,\ldots,z_d]]\longrightarrow D[[t]],\qquad
 z_j\longmapsto t^{n^{d-j}}.
\]
In two variables this is the substitution $(x,t)=(t^n,t)$.
The next lemma shows that these substitutions preserve
nondivisibility when $n$ is sufficiently large.

\begin{lemma}\label{lem:detector}
Let $a,b\in D[[z_1,\ldots,z_d]]$, where $D$ is a domain and
$a_0\ne0$. If $a\nmid b$, then
$\sigma_n(a)\nmid\sigma_n(b)$ for all sufficiently large $n$.
\end{lemma}

\begin{proof}
Write $q=b/a$ in $\Frac(D)[[z_1,\ldots,z_d]]$, and choose the
first coefficient $q_\alpha\notin D$ in lexicographic order on
$\N^d$, with the first coordinate compared first. For
$n>1+\sum_j\alpha_j$, every lexicographically later index has
weight greater than $\sum_j\alpha_j n^{d-j}$. Indeed, the increase
at the first differing coordinate exceeds the entire remaining
weight of $\alpha$. At this degree, the coefficient of
$\sigma_n(q)$ is therefore $q_\alpha$ plus a finite sum of
elements of $D$, and does not belong to $D$. Since substitution
commutes with division by a series with nonzero constant term,
$\sigma_n(a)$ does not divide $\sigma_n(b)$ in $D[[t]]$.
\end{proof}

The following lemma allows us to recover a factor over $A$ from
factors obtained after these substitutions.

\begin{lemma}\label{lem:descent}
Let $W'$ be as in Lemma~\ref{lem:division}, and let
$a,b\in W'[[z_1,\ldots,z_d]]$ satisfy $a_0=p^r$, $b_0=p^s$,
where $r,s>0$, and $ab\in A[[z_1,\ldots,z_d]]$. Suppose that
for infinitely many $n$ there is a unit $\eta_n\in W'[[t]]$ with
$(\eta_n)_0=1$ and $\eta_n\sigma_n(a)\in A[[t]]$.
Then there is a unit $\lambda\in W'[[z_1,\ldots,z_d]]$ with
$\lambda_0=1$ such that both $\lambda a$ and $\lambda^{-1}b$
belong to $A[[z_1,\ldots,z_d]]$.
\end{lemma}

\begin{proof}
We first make a one-variable observation. Suppose that
$w\in W'[[t]]$, $w_0=p^r$, and $w_j\in A$ for $j<\ell$.
If $\eta_0=1$ and $\eta w\in A[[t]]$, successive coefficient
comparison and Lemma~\ref{lem:division} give $\eta_j\in A$
for $j<\ell$. The next coefficient equation gives
\[
 w_\ell\in A+p^rW'.
\]

We now correct the coefficients of $a,b$ in lexicographic order.
Suppose that all coefficients before a nonzero index $\alpha$
already belong to $A$ in both factors. Choose one of the integers
$n$ in the hypothesis with $n>1+\sum_j\alpha_j$. Below the degree
of $\sigma_n(z^\alpha)$, and at that degree apart from
$a_\alpha$, all contributions come from earlier indices. The
observation above gives $a_\alpha\in A+p^rW'$. Multiplication
of $a$ by a suitable unit $1+\delta z^\alpha$ makes this
coefficient belong to $A$, without changing earlier coefficients.
Multiply $b$ by the inverse unit. The coefficient equation for the
unchanged product $ab$ then gives $p^r b_\alpha\in A$, and
Lemma~\ref{lem:division} gives $b_\alpha\in A$.

This construction is valid also at limit stages of the lexicographic
order. Each fixed monomial is affected only by correction monomials
dividing it, of which there are finitely many. Hence the products of
the correcting units are well defined as formal power series. After
any partial product $\lambda'$, the hypothesis remains true with
$\eta_n$ replaced by $\eta_n/\sigma_n(\lambda')$. Continuing the
construction proves the lemma.
\end{proof}

We shall also need coprimeness to persist under the two-variable
substitutions.

\begin{lemma}\label{lem:coprime}
Let $u\in W[[x,t]]$ be prime, with $u_0\ne0$, and suppose that
$u\nmid v$. Then $u(t^n,t)$ and $v(t^n,t)$ are coprime in
$W[[t]]$ for all sufficiently large $n$. The same conclusion holds
if $u,v$ are coprime products of finitely many such primes.
\end{lemma}

\begin{proof}
Otherwise, for infinitely many $n$, a common prime factor $a_n(t)$
gives a height-two prime
\[
 \mathfrak q_n=(x-t^n,a_n(t))\subset W[[x,t]]
\]
containing $u,v$. Its image in the Noetherian domain
$W[[x,t]]/(u)$ is minimal over the nonzero principal ideal
generated by $x-t^n$, and therefore has height one.

The primes $\mathfrak q_n$ are distinct. Indeed, equality for
two different indices $n,l$ would put
$t^{\min(n,l)}(1-t^{|n-l|})$, and hence $t$ and $x$, in their
common ideal. Since it also contains $u$ and $u_0\ne0$, it would
contain $p$. This is impossible for a height-two prime of the
three-dimensional regular local ring $W[[x,t]]$.

The nonzero image of $v$ would thus belong to infinitely many
height-one primes of $W[[x,t]]/(u)$. Each is minimal over $(v)$,
contradicting finiteness of the minimal primes of this ideal.
The assertion for products follows by applying the prime case to
each factor.
\end{proof}

\section{Proof of Theorem~\ref{thm:small}}

Let $A,p,W$ be as in Section~\ref{sec:prelim}. We first treat
two elementary cases of the factorization over $W[[x,t]]$.

\begin{lemma}\label{lem:two}
Let $f\in A[[x,t]]$ be irreducible, with $f_0=p^m$. Suppose that
$f\sim UV$ in $W[[x,t]]$, where $U,V$ are nonassociate primes.
Then $f$ is prime in $A[[x,t]]$.
\end{lemma}

\begin{proof}
Normalize $U_0=p^r$ and $V_0=p^{m-r}$. Suppose to the contrary
that $f\mid gh$ in $A[[x,t]]$, but $f\nmid g,h$. By
Lemma~\ref{lem:division}, these nondivisibilities remain true over
$W$. Unique factorization over $W$ therefore gives, after exchanging
$g,h$ if necessary,
\[
 U\mid g,\qquad V\mid h,\qquad U\nmid h,\qquad V\nmid g.
\]
Write $a_n=a(t^n,t)$ for a series $a(x,t)$. By
Lemma~\ref{lem:coprime}, the pairs $(U_n,h_n)$, $(V_n,g_n)$,
and $(U_n,V_n)$ are coprime over $W[[t]]$ for large $n$.

Factor $f_n$ into primes in $A[[t]]$. Every such prime divides
$g_n$ or $h_n$. In the first case, its prime factors over $W[[t]]$
all divide $U_n$; in the second, they all divide $V_n$. None of
the factors becomes a unit over $W$, since its nonunit constant
term divides $p^m$ in $A$.

Take the product, with full multiplicities, of the factors in the
first group. It belongs to $A[[t]]$ and is associated to $U_n$
over $W[[t]]$. Normalize its constant term to $p^r$, so that the
associating unit has constant term one. Lemma~\ref{lem:descent},
applied to $U$ and $f/U$, gives a proper factorization of $f$
over $A$, a contradiction.
\end{proof}

\begin{lemma}\label{lem:pure}
Let $f\in A[[x,t]]$ be irreducible, with $f_0=p^m$. If
$f\sim U^e$ in $W[[x,t]]$, where $U$ is prime and $e\geq1$,
then $f$ is prime in $A[[x,t]]$.
\end{lemma}

\begin{proof}
For $e=1$, the assertion follows from Lemma~\ref{lem:division}.
Otherwise, if $f$ is not prime, choose $g,h\in A[[x,t]]$ with
$f\mid gh$ and $f\nmid g,h$. The $U$-valuation of $g$ is then
an integer $\alpha$ with $0<\alpha<e$. Write $\alpha/e=q/N$
in lowest terms. Since $N\alpha=qe$, Lemma~\ref{lem:division}
gives $g^N=f^qk$ for some $k\in A[[x,t]]$.

Use the substitutions $(x,t)=(t^n,t)$, with subscripts as in the
preceding proof. For large $n$, Lemma~\ref{lem:coprime} makes
$U_n$ coprime to $(g/U^\alpha)_n$. In the UFD $A[[t]]$, write
\[
 a_n=\gcd(f_n,g_n),\qquad f_n=a_nb_n,\qquad g_n=a_nc_n,
 \qquad\gcd(b_n,c_n)=1.
\]
The power identity gives $b_n^q\mid a_n^{N-q}$. For a prime
$\lambda\mid U_n$ over $W[[t]]$, put $r=v_\lambda(U_n)$
and $a=v_\lambda(a_n)$. Since $a_n\mid g_n$, we have
$a\leq\alpha r$. On the other hand,
\[
 q(er-a)\leq(N-q)a
\]
gives $a\geq\alpha r$. Thus $a_n\sim U_n^\alpha$ over
$W[[t]]$. Normalize constants and apply Lemma~\ref{lem:descent}
to $U^\alpha$ and $f/U^\alpha$. This gives a proper factorization
of $f$ in $A[[x,t]]$, a contradiction.
\end{proof}

The next lemma applies when all factors over the completed
coefficient ring have constant term $p$.

\begin{lemma}\label{cor:split}
Let $W'$ be a complete DVR extension of $W$ in which $p$ remains
a uniformizer. Let $f\in A[[x,t]]$ be irreducible, with $f_0=p^m$.
If every prime factor of $f$ in $W'[[x,t]]$ has constant term
associated to $p$, then $f$ is prime in $A[[x,t]]$.
\end{lemma}

\begin{proof}
Write $f=\varepsilon\prod_{i=1}^m F_i$ over $W'$, listing factors
with multiplicity and normalizing $(F_i)_0=p$, $\varepsilon_0=1$.
Suppose that $f\mid gh$ in $A[[x,t]]$, but $f\nmid g,h$.
By Lemmas~\ref{lem:division} and~\ref{lem:detector}, for all
large $n$ we have $f_n\nmid g_n,h_n$. Each $(F_i)_n$ is prime
in $W'[[t]]$, since its constant term is $p$ and it has order
one in that regular local ring.

Since $A[[t]]$ is a UFD, $f_n$ has a proper divisor $a_n$.
Normalize $(a_n)_0=p^{r_n}$, where $0<r_n<m$. Over $W'[[t]]$,
this divisor is associated to a proper nonempty subproduct of the
listed $(F_i)_n$. There are only finitely many choices of the
indices, so, on an infinite subsequence, the same product
$D=\prod_{i\in I}F_i$ occurs. The associating units have constant
term one. Lemma~\ref{lem:descent}, applied to $D$ and $f/D$,
now contradicts irreducibility of $f$ over $A$.
\end{proof}

\begin{proof}[Proof of Theorem~\ref{thm:small} for $m\leq3$]
Normalize $f_0=p^m$ and list its prime factors over $W[[x,t]]$
with multiplicity. Their positive constant-term $p$-valuations form
a partition of $m$. For $m\leq3$, either there is only one factor,
all factors have constant term associated to $p$, or the partition
is $(2,1)$. The first case follows from Lemma~\ref{lem:division},
the second from Lemma~\ref{cor:split}, and the third from
Lemma~\ref{lem:two}.
\end{proof}

We now prepare for the case $m=4$. Recall that an irreducible series
$h\in W[[t]]$ is irreducible modulo $t^j$ for some $j$: no product
of two nonunits is congruent to $h$ modulo $t^j$
\cite[Theorem~6.4]{JP}. This finite irreducibility theorem holds
over every Krull domain; its characteristic-zero UFD case was proved
by Bayart \cite[\S IV, Th\'eor\`eme~(T)]{Bayart81}.
We shall use the following consequence to approximate coefficients.

\begin{lemma}\label{lem:stability}
Let $h\in W[[t]]$ satisfy $h_0=p^m$, and suppose that $h$ is
irreducible modulo $t^j$, where $j\geq2$. Every series congruent
to $h$ modulo $(p^{mj},t^j)$ is irreducible in $W[[t]]$.
\end{lemma}

\begin{proof}
Write $p^m=h-th_1$. Taking $j$th powers gives
\[
 p^{mj}=hQ+t^jR,\qquad Q_0=p^{m(j-1)}.
\]
If $h'=h+p^{mj}c+t^jd$, then $h'\equiv h(1+cQ)\pmod{t^j}$.
The constant term of $1+cQ$ is $1+c_0p^{m(j-1)}$, which is a
unit of $W$. Hence $h'$ is irreducible modulo $t^j$, and therefore
irreducible.
\end{proof}

For a quadratic factor, we shall need a substitution preserving
irreducibility after extending the residue field.

\begin{lemma}\label{lem:quadratic-odd}
Suppose that $W=k[[p]]$, where $k$ is algebraically closed of
characteristic zero. If $\Phi\in W[[x,t]]$ is prime and
$\Phi_0=p^2$, then there are $\alpha,\beta\in tk[t]$ such that
$\Phi(\alpha(t),\beta(t))$ is irreducible in $W[[t]]$.
\end{lemma}

\begin{proof}
By Weierstrass preparation in $p$,
\[
 \Phi\sim p^2+a(x,t)p+b(x,t),\qquad
 a,b\in(x,t)k[[x,t]].
\]
The discriminant $\Delta=a^2-4b$ is not a square in
$\Frac(k[[x,t]])$. Indeed, a square root in the fraction field
would belong to $k[[x,t]]$ by integral closedness, and would give
a proper factorization of the distinguished quadratic.

Every unit of $k[[x,t]]$ is a square. Thus some irreducible factor
$\pi$ of $\Delta$ occurs with odd multiplicity $e$. Write
$\Delta=\pi^eH$, where $\pi\nmid H$. Normalize the plane branch
$\pi=0$, obtaining a parametrization
\[
 x=\alpha_0(z),\qquad t=\beta_0(z).
\]
Its normalization is a finite complete DVR with residue field $k$,
and hence is $k[[z]]$. The kernel of the parametrization is $(\pi)$,
so $H(\alpha_0,\beta_0)\ne0$.

At least one partial derivative of $\pi$ is nonzero on this branch.
To see this, make a linear change of coordinates and prepare $\pi$
as a unit times a monic polynomial in $t$. Its $t$-derivative is
nonzero in characteristic zero and has smaller degree, so it cannot
be divisible by that polynomial. Nonvanishing of the gradient on
the branch is unchanged by the coordinate change.

Suppose, for example, that $\pi_t(\alpha_0,\beta_0)$ has order $c$.
For $L>c$, expansion in the second variable gives
\[
 \ord_z\pi(\alpha_0,\beta_0+z^L)=L+c,
\]
since every higher-order term has order at least $2L$. For large
$L$, the order of $H$ on this perturbed parametrization remains
a fixed integer $d$. Consequently
\[
 \ord_z\Delta(\alpha_0,\beta_0+z^L)=e(L+c)+d.
\]
Choose the parity of $L$ to make this number odd. The substituted
discriminant is then not a square in $k((z))$, so the substituted
quadratic is irreducible in $k[[p,z]]$, and the same is true of
the corresponding substitution in $\Phi$.
Truncating $\alpha_0$ and $\beta_0+z^L$ to polynomials preserves
the first nonzero discriminant coefficient and gives the required
$\alpha,\beta$.
\end{proof}

We next choose a finite coefficient extension over which to take
the norm.

\begin{lemma}\label{lem:coefficientmodel}
Suppose that $W=k[[p]]$, and let $\ell/k$ be a finite Galois
extension of degree $d$. There is $\eta\in\ell[[p]]$ satisfying
a monic polynomial of degree $d$ over $A$, such that
$A[\eta]$ is free of rank $d$ over $A$,
\[
 \ell[[p]]=W[\eta],
\]
and the conjugates of $\eta$ have pairwise distinct residues.
Moreover, any finite family of polynomials in $t\ell[[p]][t]$
can be approximated to any prescribed $p$-adic precision by
elements of $tA_{(p)}[\eta][t]$. A subsequent substitution
$t\mapsto st$, with $s\in A\setminus(p)$, puts all their
coefficients in $A[\eta]$.
\end{lemma}

\begin{proof}
Choose a primitive element $\gamma$ of $\ell/k$, and lift its
monic minimal polynomial to a polynomial $q$ over $A_{(p)}$.
Hensel's lemma in $\ell[[p]]$ gives a root $\eta_0$ of $q$
reducing to $\gamma$. The elements
$1,\eta_0,\ldots,\eta_0^{d-1}$ form a $W$-basis of $\ell[[p]]$:
their residues form a $k$-basis of $\ell$, so their change-of-basis
matrix has unit determinant.

Write the coefficients of $q$ with a common denominator
$s_0\in A\setminus(p)$. For sufficiently large $N$,
\[
 q_1(z)=s_0^{Nd}q(z/s_0^N)
\]
is monic of degree $d$ over $A$. Put $\eta=s_0^N\eta_0$.
Its powers remain a $W$-basis, and are therefore linearly independent
over $A$. Thus $A[\eta]$ is free of rank $d$ and
$W[\eta]=\ell[[p]]$. The conjugate residues remain distinct,
since $s_0$ is a unit of $W$.

Expand the finitely many coefficients to be approximated in this
basis and use the density of $A_{(p)}$ in $W$. The resulting
coefficients have finitely many denominators in $A\setminus(p)$.
A common multiple $s$ of these denominators clears them after
$t\mapsto st$, since all terms have positive $t$-degree. This
substitution is an automorphism over $\ell[[p]]$.
\end{proof}

\Needspace{22\baselineskip}
The next lemma is the norm argument. The finite coefficient
extension in its statement is not assumed to have a UFD power
series ring.

\begin{lemma}\label{thm:norm}
Let $f\in A[[x,t]]$ be irreducible, with $f_0=p^{s+r}$, and
suppose that
\[
 f\sim U_1\cdots U_s\Phi\quad\text{in }W[[x,t]],\qquad
 (U_i)_0=p,\quad\Phi_0=p^r,\quad r\geq2,
\]
where the displayed factors are prime. The $U_i$ may repeat, and
$\Phi$ occurs once. Suppose that there is a finite unramified
Galois extension $W'/W$ and $\eta\in W'$ such that:
\begin{enumerate}
\item $W'=W[\eta]$, and $\eta$ satisfies a monic polynomial over
      $A$ of degree $d=[\Frac(W'):\Frac(W)]$;
\item the conjugates of $\eta$ have pairwise distinct residues;
\item there are $\alpha,\beta\in tA[\eta][t]$ such that
      $\Phi(\alpha(t),\beta(t))$ is irreducible in $W'[[t]]$.
\end{enumerate}
Then $f$ is prime in $A[[x,t]]$.
\end{lemma}

\begin{proof}
Put $A'=A[\eta]$. The monic polynomial in (i) is the minimal
polynomial of $\eta$ over $\Frac(W)$, so $A'$ is free of rank
$d$ over $A$. Hence $A'[[u,v,a,b,z]]$ is finite free of rank $d$
over $A[[u,v,a,b,z]]$, and its determinant norm takes values in
the latter ring.

Suppose to the contrary that $f\mid gh$ in $A[[x,t]]$, but
$f\nmid g,h$. Put
\[
 f^\dagger=f\bigl(u+\alpha(z)+\eta a,\,
                       v+\beta(z)+\eta b\bigr),
\]
and define $g^\dagger,h^\dagger$ in the same way. Denote their
norms by $n_f,n_g,n_h$. Then
\[
 n_f\mid n_gn_h,\qquad
 n_f\big|_{a=b=z=0}=f(u,v)^d,
\]
with analogous identities for $g,h$. Even over $W'$, $n_f$ divides
neither $n_g$ nor $n_h$. Indeed, a divisibility would give
$f^d\mid g^d$ or $f^d\mid h^d$ after setting $a=b=z=0$.
Unique factorization over $W'$ and Lemma~\ref{lem:division}
would then give $f\mid g$ or $f\mid h$ over $A$.

Write $\Gamma=\operatorname{Gal}(W'/W)$, and use subscripts
$\sigma$ for the conjugates of $\eta,\alpha,\beta$. Over $W'$,
the factors of $n_f$ are, up to a unit,
\begin{align*}
 U_{i,\sigma}&=U_i\bigl(u+\alpha_\sigma(z)+\eta_\sigma a,\,
                            v+\beta_\sigma(z)+\eta_\sigma b\bigr),\\
 \Phi_\sigma&=\Phi\bigl(u+\alpha_\sigma(z)+\eta_\sigma a,\,
                            v+\beta_\sigma(z)+\eta_\sigma b\bigr).
\end{align*}
The $U_i$ remain prime over $W'$, since their constant terms are
$p$. Also $\Phi$ remains prime: a proper factorization over $W'$
would give a proper factorization after the substitution in (iii),
since both factor constants would have positive $p$-valuation.
The displayed translations are automorphisms after adjoining
$a,b,z$, so all the $U_{i,\sigma}$ and $\Phi_\sigma$ are prime.

The factors $\Phi_\sigma$ are pairwise nonassociate. Otherwise,
set $z=0$ and undo one of the two translations to obtain
\[
 \Phi(u+\delta a,v+\delta b)\sim\Phi(u,v),\qquad
 \delta=\eta_\sigma-\eta_\tau\in(W')^\times.
\]
Setting $a=-u/\delta$ and $b=-v/\delta$ gives
$p^r\sim\Phi(u,v)$, contradicting primality of $\Phi$ and
$r\geq2$. The group $\Gamma$ acts transitively on these $d$ factors.

Now substitute
\[
 (u,v,a,b,z)=(t^{n^4},t^{n^3},t^{n^2},t^n,t),
\]
and denote the resulting series by a subscript $n$. Each factor
is congruent modulo $t^n$ to its corresponding substitution from
(iii), or its conjugate. For $U_i$ the resulting constant term is
$p$, so those substitutions are irreducible as well. Finite
irreducibility makes all factor substitutions prime in $W'[[t]]$
for large $n$. By Lemma~\ref{lem:detector}, the distinct
$(\Phi_\sigma)_n$ remain nonassociate. None is associated to
a $(U_{i,\tau})_n$, since their constant-term valuations are
$r$ and $1$.

The same lemma, applied to $n_f,n_g,n_h$, gives
\[
 (n_f)_n\mid(n_g)_n(n_h)_n,\qquad
 (n_f)_n\nmid(n_g)_n,\qquad (n_f)_n\nmid(n_h)_n
\]
for large $n$. Hence $(n_f)_n$ is reducible in $A[[t]]$.
Every prime in its factorization over $A[[t]]$ remains a nonunit
over $W'$, since its constant term is a unit of $A$ times a
positive power of $p$.

The factorization over $W'[[t]]$ of a series in $A[[t]]$ is
$\Gamma$-invariant. Thus a prime factor over $A[[t]]$ containing
one $(\Phi_\sigma)_n$ must contain all of them. Since each occurs
once, exactly one prime factor contains them, and it occurs once.
Reducibility gives another prime factor $c_n\in A[[t]]$, whose
factors over $W'$ all come from the $(U_{i,\sigma})_n$.

There are finitely many subproducts of the indexed list
$U_{i,\sigma}$. On an infinite subsequence, $c_n$ is associated
to the substitution of a fixed nonempty product $D$ of these
factors. If $D_0=p^e$, normalize $(c_n)_0=p^e$. The associating
units then have constant term one. Apply Lemma~\ref{lem:descent}
in five variables to $n_f=D(n_f/D)$. It gives a divisor
$c\in A[[u,v,a,b,z]]$ of $n_f$, associated to $D$ over $W'$.

Set $a=b=z=0$. The resulting nonunit $c^*\in A[[u,v]]$ satisfies
\[
 c^*\mid f(u,v)^d,\qquad
 c^*\sim\prod_{i=1}^s U_i^{e_i}\quad\text{over }W'[[u,v]],
 \qquad\sum_i e_i=e>0.
\]
Factor $c^*$ into irreducibles over $A$, which is possible by
induction on the $p$-valuation of its constant term. Each irreducible
has only factors with constant term associated to $p$ over $W'$,
and is therefore prime by Lemma~\ref{cor:split}. Such a prime
divides $f^d$, hence divides $f$. It must be associated to the
irreducible $f$, which is impossible because its factorization over
$W'$ does not contain $\Phi$. This contradiction proves the lemma.
\end{proof}

\begin{proof}[Completion of the proof of Theorem~\ref{thm:small}]
It remains to consider $m=4$ and $A\supseteq\Q$. Normalize
$f_0=p^4$. By the Cohen structure theorem,
$W\cong k[[p]]$, using the specified uniformizer $p$
\cite[Tag 0323]{Stacks}. The partitions $(4)$, $(3,1)$, and
$(2,2)$ of the constant-term valuation are covered by
Lemmas~\ref{lem:two} and~\ref{lem:pure}, and $(1,1,1,1)$ by
Lemma~\ref{cor:split}. We are left with
\[
 f\sim U_1U_2\Phi,\qquad (U_1)_0=(U_2)_0=p,\quad\Phi_0=p^2,
\]
where $U_1,U_2$ may coincide and $\Phi$ occurs once.

Extend coefficients to $\overline{k}[[p]]$. If $\Phi$ becomes
reducible, both its factors have constant term associated to $p$,
so Lemma~\ref{cor:split} applies. Otherwise,
Lemma~\ref{lem:quadratic-odd} gives a polynomial substitution over
$\overline{k}$ for which $\Phi$ remains irreducible. The finitely
many coefficients lie in a finite Galois extension $\ell/k$.
Choose $A[\eta]$ as in Lemma~\ref{lem:coefficientmodel}, and
$j\geq2$ such that this substitution of $\Phi$ is irreducible
modulo $t^j$.

Approximate its coefficients by elements of $A_{(p)}[\eta]$
modulo $p^{2j}$. Lemma~\ref{lem:stability} preserves
irreducibility. After a substitution $t\mapsto st$,
$s\in A\setminus(p)$, both polynomials belong to
$tA[\eta][t]$ and the substitution of $\Phi$ remains irreducible
over $\ell[[p,t]]$. Lemma~\ref{thm:norm} now applies.
\end{proof}

\section{Proof of Theorem~\ref{thm:reduction}}
\label{sec:reduction}

Put $D=A[u,v]$ and $B=D[[t]]$, and suppose that $B$ is a UFD.
For each integer $j$, let $D_j$ be the polynomials homogeneous
of total degree $j$ in $u,v$, with $D_j=0$ for $j<0$. For an
integer $w$, put
\[
 \Gamma_w=\left\{\sum_{n\geq0}b_nt^n\in B:
                          b_n\in D_{n+w}\right\}.
\]
Thus $u,v$ have weight one, $t$ has weight minus one, and
multiplication adds weights. The substitution $x\mapsto ut$,
$t\mapsto vt$ gives an isomorphism
\begin{equation}\label{red:eq:embedding}
 A[[x,t]]\longrightarrow\Gamma_0,\qquad
 \sum_{i,j\geq0}c_{ij}x^it^j\longmapsto
 \sum_{n\geq0}\left(\sum_{i+j=n}c_{ij}u^iv^j\right)t^n.
\end{equation}

For $b=\sum_{n\geq0}b_nt^n\in B$, denote its component of
weight $w$ by
\[
 b^{[w]}=\sum_{n\geq0}(b_n)_{n+w}t^n,
\]
where $(b_n)_j$ is the degree-$j$ homogeneous component of $b_n$.
Only finitely many weights occur modulo $t^N$ for any fixed $N$,
so components may be compared in products of formal power series.

\begin{lemma}\label{red:lem:cancellation}
If $0\ne h\in\Gamma_d$ and $ha\in\Gamma_e$ for $a\in B$,
then $a\in\Gamma_{e-d}$.
\end{lemma}

\begin{proof}
The component of $ha$ of weight $d+w$ is $ha^{[w]}$.
For $w\ne e-d$ it is zero. Since $B$ is a domain, all such
$a^{[w]}$ are zero.
\end{proof}

We shall show that each prime factor of an element of $\Gamma_0$,
apart from $t$, has a homogeneous associate.

\begin{lemma}\label{red:lem:components}
Let $0\ne F\in\Gamma_0$, and let $q$ be a prime factor of $F$
in $B$, with $q_0\ne0$. Then $q$ divides $q^{[w]}$ in $B$
for every integer $w$.
\end{lemma}

\begin{proof}
Introduce a formal variable $z$ and define
\[
 \Delta:B\longrightarrow B[[z]],\qquad
 u\longmapsto(1+z)u,\quad v\longmapsto(1+z)v,\quad
 t\longmapsto(1+z)^{-1}t.
\]
This is a well-defined homomorphism fixing $A$, and
$\Delta(F)=F$. Write $F=q^eH$, where $q\nmid H$. Reducing
$\Delta(F)=\Delta(q)^e\Delta(H)$ modulo $q$ gives zero in
the domain $(B/(q))[[z]]$. The constant $z$-coefficient of
$\Delta(H)$ is $H\notin(q)$, so
\begin{equation}\label{red:eq:invariant-principal}
 \Delta(q)\in qB[[z]].
\end{equation}

Fix $N\geq1$, and put $C_N=B/(q,t^N)$. Only finitely many
weights of $q$ remain modulo $t^N$, so
\eqref{red:eq:invariant-principal} gives
\[
 0=\sum_w\overline{q^{[w]}}(1+z)^w\quad\text{in }C_N[[z]].
\]
Distinct integral powers of $1+z$ are linearly independent over
any coefficient ring: multiply a finite relation by a large power
of $1+z$ and apply the polynomial-ring isomorphism sending a
variable to $1+z$. Thus $q^{[w]}\in(q,t^N)$ for every $w,N$.

The ideal $(q)$ is $t$-adically closed, since $q_0\ne0$.
Indeed, if $b\in\bigcap_N(q,t^N)$, choose $a_N\in B$ with
$b-qa_N\in t^NB$. For $M\geq N$,
$q(a_M-a_N)\in t^NB$. Since $q$ has $t$-order zero and
$D$ is a domain, $a_M-a_N\in t^NB$. The sequence $a_N$
converges in $B$, and its limit $a$ satisfies $b=qa$.
It follows that $q^{[w]}\in qB$.
\end{proof}

\begin{proof}[Proof of Theorem~\ref{thm:reduction}]
Let $0\ne F\in\Gamma_0$. Factor it in $B$ as
\begin{equation}\label{red:eq:ambient}
 F=\varepsilon t^m\prod_{i=1}^k q_i^{e_i},
\end{equation}
where $\varepsilon\in B^\times$ and the $q_i$ are prime and
not associated to $t$. Then $(q_i)_0\ne0$. The coefficient of
$t^m$ is
\[
 F_m=\varepsilon_0\prod_i(q_i)_0^{e_i}\in D_m,
 \qquad\varepsilon_0\in D^\times=A^\times.
\]
A product of nonzero polynomials over a domain is homogeneous
only if each factor is homogeneous, by comparison of the least
and greatest total degrees. Thus $(q_i)_0$ is homogeneous of
some degree $d_i\geq0$, and
\begin{equation}\label{red:eq:weight-balance}
 m=\sum_i e_id_i.
\end{equation}

By Lemma~\ref{red:lem:components}, the homogeneous component
$h_i=q_i^{[d_i]}$ is divisible by $q_i$, say $h_i=q_iv_i$.
Since $(h_i)_0=(q_i)_0\ne0$, we have $(v_i)_0=1$. Hence
$h_i$ is a prime associate of $q_i$ in $B$.
Put $P_i=t^{d_i}h_i\in\Gamma_0$. Lemma~\ref{red:lem:cancellation}
and $\Gamma_{-d_i}=t^{d_i}\Gamma_0$ give
\begin{equation}\label{red:eq:contraction}
 h_iB\cap\Gamma_0=P_i\Gamma_0.
\end{equation}
Indeed, if $h_iH\in\Gamma_0$ for $H\in B$, then
$H\in\Gamma_{-d_i}=t^{d_i}\Gamma_0$. The ideal on the left
of \eqref{red:eq:contraction} is prime, so $P_i$ is prime
in $\Gamma_0$.

Replacing $q_i$ by $h_i$ in \eqref{red:eq:ambient} changes
only the unit. By \eqref{red:eq:weight-balance}, the product
of $t^m$ and the $h_i^{e_i}$ has weight zero. Thus the remaining
unit belongs to $\Gamma_0$ by Lemma~\ref{red:lem:cancellation},
and so does its inverse. Thus $F=\varepsilon'\prod_i P_i^{e_i}$
is a factorization into primes in $\Gamma_0$. Therefore $\Gamma_0$, and hence
$A[[x,t]]$, is a UFD.

Finally, suppose that \eqref{eq:pep} holds for every domain $B$,
and that $A[[t]]$ is a UFD. Apply \eqref{eq:pep} first to $A$
and then to $A[u]$. It follows that $A[u,v][[t]]$ is a UFD,
so the first assertion gives that $A[[x,t]]$ is a UFD.
\end{proof}

\subsection*{Declaration of AI use}
This paper was written with AI assistance, model ChatGPT 5.6-Sol. The Proof of Theorem \ref{thm:small} was obtained by the author for $m \leq 3$; The extension for $m = 4$, and the reduction given by Theorem \ref{thm:reduction}, were obtained by the AI model.


\begin{thebibliography}{99}
\raggedright
\bibitem{Bayart73}
M. Bayart, \emph{S\'eries formelles sur un anneau factoriel},
C. R. Acad. Sci. Paris S\'er. A \textbf{277} (1973), 449--450.

\bibitem{Bayart80}
M. Bayart, \emph{Factorialit\'e et s\'eries formelles irr\'eductibles I},
in S\'eminaire d'Alg\`ebre Paul Dubreil et Marie-Paule Malliavin,
Lecture Notes in Mathematics \textbf{795},
Springer, Berlin, 1980, 402--433.
\href{https://doi.org/10.1007/BFb0090130}{doi:10.1007/BFb0090130}.

\bibitem{Bayart81}
M. Bayart, \emph{Factorialit\'e et s\'eries formelles irr\'eductibles II},
in S\'eminaire d'Alg\`ebre Paul Dubreil et Marie-Paule Malliavin,
Paris 1980, Lecture Notes in Mathematics \textbf{867},
Springer, Berlin, 1981, 174--213.
\href{https://doi.org/10.1007/BFb0090387}{doi:10.1007/BFb0090387}.

\bibitem{Buchsbaum}
D. A. Buchsbaum, \emph{Some remarks on factorization in power series rings},
J. Math. Mech. \textbf{10} (1961), 749--753.
\href{https://doi.org/10.1512/iumj.1961.10.10052}{doi:10.1512/iumj.1961.10.10052}.

\bibitem{Claborn}
L. Claborn, \emph{Note generalizing a result of Samuel's},
Pacific J. Math. \textbf{15} (1965), 805--808.
\url{https://msp.org/pjm/1965/15-3/pjm-v15-n3-p07-p.pdf}.

\bibitem{DanilovSamuel}
V. I. Danilov, \emph{On a conjecture of Samuel},
Math. USSR-Sb. \textbf{10} (1970), 127--137.
\href{https://doi.org/10.1070/SM1970v010n01ABEH001590}{doi:10.1070/SM1970v010n01ABEH001590}.

\bibitem{DanilovDGC}
V. I. Danilov, \emph{Rings with a discrete group of divisor classes},
Math. USSR-Sb. \textbf{12} (1970), 368--386.
\href{https://doi.org/10.1070/SM1970v012n03ABEH000926}{doi:10.1070/SM1970v012n03ABEH000926}.

\bibitem{JP}
A. Jones and E. Paran, \emph{On Landweber's unique factorization
problem}, arXiv:2607.03475v2 (2026).
\url{https://arxiv.org/abs/2607.03475}.

\bibitem{Krull}
W. Krull, \emph{Beitr\"age zur Arithmetik kommutativer
Integrit\"atsbereiche. V. Potenzreihenringe},
Math. Z. \textbf{43} (1938), 768--782.
\href{https://doi.org/10.1007/BF01181119}{doi:10.1007/BF01181119}.

\bibitem{Lipman}
J. Lipman, \emph{Rings with discrete divisor class group:
theorem of Danilov--Samuel},
Amer. J. Math. \textbf{101} (1979), 203--211.
\href{https://doi.org/10.2307/2373945}{doi:10.2307/2373945}.

\bibitem{Salmon}
P. Salmon, \emph{Sur les s\'eries formelles restreintes},
Bull. Soc. Math. France \textbf{92} (1964), 385--410.
\href{https://doi.org/10.24033/bsmf.1613}{doi:10.24033/bsmf.1613}.

\bibitem{Samuel}
P. Samuel, \emph{On unique factorization domains},
Illinois J. Math. \textbf{5} (1961), 1--17.
\href{https://doi.org/10.1215/ijm/1255629643}{doi:10.1215/ijm/1255629643}.

\bibitem{Stacks}
The Stacks Project,
\url{https://stacks.math.columbia.edu/}, 2026;
Tag \href{https://stacks.math.columbia.edu/tag/0323}{0323}
(Cohen structure).

\bibitem{Tran}
V.-H. Tran, D. V. Nguyen, Q. X. Nguyen, T. N. Vo, and T. M. Nguyen,
\emph{The Noetherian case of Bayart's power-series question},
arXiv:2608.12642 (2026).
\url{https://arxiv.org/abs/2608.12642}.
\end{thebibliography}
\end{document}